\documentclass{amsart}
\usepackage{amssymb}

\theoremstyle{plain}

\numberwithin{equation}{section}

\begin{document}

\date{October 20, 2001}

\noindent {\Large \bf Op\'erades Lie-admissibles}

\medskip

\noindent \textbf{Elisabeth REMM}

\bigskip

\noindent \textbf{Laboratoire de Math\'{e}matiques et Applications, 4, rue
des Fr\`{e}res Lumi\`{e}re, 68093 Mulhouse cedex, France}

\noindent \textbf{E-mail : E.Remm@uha.fr}

\bigskip

\bigskip

\textbf{Abstract.} The aim of this paper is to present remarkable classes of 
Lie-admissible algebras containing in particular the associative algebras, the Vinberg and pre-Lie algebras. We determine the 
associated operads and their dual operads. 
\smallskip

\textbf{R\'{e}sum\'{e}. }Le but de cette note est de pr\'{e}senter des
classes remarquables d'alg\`{e}bres Lie-admissibles qui contiennent entre
autres les alg\`{e}bres associatives, de Vinberg et pr\'{e}-Lie
et de d\'{e}terminer leurs op\'{e}rades associ\'{e}es et les op\'{e}rades
duales.

\bigskip

\noindent 1. {\bf Alg\`{e}bres Lie-Admissibles.}

\smallskip

Soient $(A,\mu )$ une $\Bbb{K}$-alg\`ebre sur un corps commutatif de
caract\'{e}ristique 0. Notons $a_{\mu }:A^{\otimes 3}\rightarrow A$
l'associateur de la loi $\mu $ :  \[
a_{\mu }\left( X_1,X_2,X_3\right) =\mu \left( \mu \left( X_1,X_2\right)
,X_3\right) -\mu \left( X_1,\mu \left( X_2,X_3\right) \right).
\]
Soit $\sum_{n}$ le groupe sym\'{e}trique d'ordre $n$. Pour tout $\sigma \in
\sum_{3}$, on pose
\[
\sigma \left( X_{1},X_{2},X_{3}\right) = \left( X_{\sigma ^{-1}\left( 1\right) },X_{\sigma ^{-1}\left( 2\right)
},X_{\sigma ^{-1}\left( 3\right) }\right).
\]

\noindent DEFINITION 1 : L'alg\`{e}bre $\mathcal{A=}\left( A,\mu \right) $ est
dite Lie-admissible (cf [1]) si sa loi $\mu $ v\'erifie
\[
\sum_{\sigma \in \sum_{3}}\left( -1\right) ^{\varepsilon \left( \sigma
\right) }a_{\mu } \circ \sigma=0.
\]

\medskip

\noindent Soit $\mu $ une loi d'alg\`{e}bre Lie-admissible; alors 
\[
\left[ X,Y\right] _{\mu }:=\mu \left( X,Y\right) -\mu \left( Y,X\right)
\]
est une loi d'alg\`{e}bre de Lie. Si $\mathcal{A=}\left( A,\mu \right) $ est
Lie-admissible, on notera $\mathcal{A}_{L}=(A,[,]_\mu) $ l'alg\`{e}bre de Lie
ainsi d\'efinie.

\medskip

\noindent {\bf Remarque}. Pour toute alg\`{e}bre de Lie $\frak{g}$, il existe
une alg\`{e}bre Lie-admissible $\mathcal{A=} \left( A,\mu \right) $ telle que
$\frak{g}=\mathcal{A}_{L}.$ En effet
$\mu (X,Y)= \frac{1}{2} [X,Y]$ est une loi Lie-admissible v\'erifiant la condition
demand\'ee. Notons \'egalement que dans le cas de la dimension finie, la d\'eriv\'ee covariante ([6]) d'une connection 
de Levi Civita associ\'ee \`a une forme quadratique d\'efinie positive d\'etermine \'egalement une
loi Lie-admissible. 

\bigskip

\noindent 2. {\bf Classes d'alg\`{e}bres Lie-admissibles.}

\smallskip

\noindent 2.1. Alg\`{e}bres $G_{i}$-associatives.

\smallskip

Notons $G_{1}=\left\{
Id\right\} $, $G_{2}=\left\{ Id,\tau _{12}\right\} $, $G_{3}=\left\{ Id,\tau
_{23}\right\} $, $G_{4}=\left\{ Id,\tau _{13}\right\} $, $G_{5}=\left\{
Id,3-cycles\right\}=\mathcal{A}_{3}$, et $G_{6}=\sum_{3}$ les
diff\'erents sous-groupes de $\sum_{3}.$

\medskip

\noindent DEFINITION\ 2 : Une $\Bbb{K}$-alg\`ebre 
$\mathcal{A=}\left( A,\mu \right) $ est dite $G_{i}$-associative si sa loi 
$\mu $ v\'{e}rifie \[
\sum_{\sigma \in G_{i}}\left( -1\right) ^{\varepsilon \left( \sigma \right)
} a_{\mu } \circ \sigma=0. \qquad \left( \ast \right)
\]

\medskip

\noindent Notons qu'une alg\`{e}bre $\left( A,\mu \right) $ v\'{e}rifiant
$\left( \ast \right) $ est n\'{e}cessairement Lie-admissible.

\medskip

\noindent 2.2. Description.

\smallskip

\noindent a) $G_{1}=\left\{ Id\right\} $. Les alg\`{e}bres correspondantes
sont les alg\`{e}bres associatives.

\noindent b) $G_{2}=\left\{ Id,\tau _{12}\right\} $. La relation $\left( \ast
\right) $ s'\'{e}crit \[
\mu \left( \mu \left( X_{1},X_{2}\right) ,X_{3}\right) -\mu \left( X_{1},\mu
\left( X_{2},X_{3}\right) \right) -\mu \left( \mu \left( X_{2},X_{1}\right)
,X_{3}\right) +\mu \left( X_{2},\mu \left( X_{1},X_{3}\right) \right) =0
\]
Les alg\`{e}bres $G_{2}$-associatives  sont les alg\`ebres de Vinberg (appel\'ees aussi alg\`ebres sym\'{e}triques-gauche) [10]
. Si $\mathcal{A}$ est de dimension
finie, l'alg\`{e}bre de Lie associ\'{e}e $\mathcal{A}_{L}$ est munie d'une
structure affine.

\noindent c) $G_{3}=\left\{ Id,\tau _{23}\right\} $ Les alg\`{e}bres
correspondantes sont les alg\`{e}bres pr\'{e}-Lie (voir [5]). Ces alg\`ebres sont aussi appel\'ees sym\'etriques-droite. Notons que si la loi 
$x.y$ est pr\'{e}-Lie, alors la loi $x \odot y = y.x$ est de Vinberg. 

\noindent Les alg\`{e}bres $G_{4}$ et $G_{5}$-associatives ne semblent pas
avoir fait l'objet d'\'{e}tudes particuli\`{e}res. Notons que si $\mu $ est
$G_{5}$-associative, elle v\'{e}rifie :
\begin{eqnarray*}
&&\mu \left( \mu \left( X_{1},X_{2}\right) ,X_{3}\right) -\mu \left(
X_{1},\mu \left( X_{2},X_{3}\right) \right) +\mu \left( \mu \left(
X_{2},X_{3}\right) ,X_{1}\right)  \\
&&-\mu \left( X_{2},\mu \left(
X_{3},X_{1}\right) \right)  +\mu \left( \mu \left( X_{3},X_{1}\right) ,X_{2}\right) -\mu \left(
X_{3},\mu \left( X_{1},X_{2}\right) \right) = 0.  
\end{eqnarray*}
Si l'on suppose de plus que $\mu $ est anticommutative$\left( \left[
,\right] _{\mu }=2\mu \right) $ la relation ci-dessus est la condition de
Jacobi de la loi d'alg\`{e}bre de Lie $\mu .$ Dans l'optique d'une
des- \-cription d'une version ''non-commutative'' des alg\`{e}bres de Lie
donn\'{e}e par exemple par les alg\`ebres de Leibniz [8], les alg\`{e}bres
$G_{5}$-associatives peuvent \^{e}tre \'{e}galement consid\'{e}r\'{e}es comme
des candidats naturels.

\bigskip

\noindent 3. {\bf Op\'{e}rades associ\'{e}es aux alg\`{e}bres
$G_{i}$-associatives}

\smallskip

Soit $\Bbb{K}\left[ \sum_{n}\right]$ la $\Bbb{K}$-alg\`ebre du groupe sym\'etrique $\sum_{n}$. Une op\'erade $\mathcal{P}$ 
est d\'efinie par une suite d'espaces vectoriels sur $\Bbb{K}$, $\mathcal{P}(n)$, $n\geq1$ telle que $\mathcal{P}(n)$ soit un module sur 
$\Bbb{K}\left[ \sum_{n}\right]$ et par des applications de composition
\[
\circ _{i}:\mathcal{P}(n)\otimes \mathcal{P}(m)\rightarrow \mathcal{P}%
(n+m-1)\qquad i=1,...,n 
\]
satisfaisant des propri\'{e}t\'{e}s ''associatives'', les axiomes de May [9].

Tout $\Bbb{K}\left[ \sum_{n}\right] $-module $E$ engendre une op\'erade libre not\'ee $\mathcal{F(}E)$ ([5])
v\'erifiant $\mathcal{F(}E)(1)=\Bbb{K}$, $\mathcal{F(}E)(2)=E$.
En particulier si $E= \Bbb{K}\left[ \sum_{2}\right] $, le module libre $\mathcal{F} \left( E\right) (n)$ admet comme base 
les ''produits parenth\'{e}s\'{e}s'' de $n$ variables index\'{e}es par
$\left\{ 1,2,...,n\right\}$ .

Soit $E$ un $\Bbb{K}\left[ \sum_{2}\right]$-module et $R$ un $\Bbb{K}\left[ \sum_{3}\right]$-sous-module de $\mathcal{F(}E)(3)$. On note 
$\mathcal{R}$ l'id\'eal engendr\'e par $R$, c'est-\`a-dire l'intersection de tous les id\'eaux $\mathcal{I}$ de $\mathcal{F}(E)$ tels que 
$\mathcal{I}(1)=\left\{ 0 \right\}$, $\mathcal{I}(2)=\left\{0 \right\}$ et $\mathcal{I}(3)=R$.

On appelle op\'erade binaire quadratique engendr\'ee par $E$ et par les relations $R$ l'op\'erade $\mathcal{P}(\Bbb{K},E,R)$, not\'ee \'egalement 
$\mathcal{F}(E)/\mathcal{R}$, d\'efinie par :
\[
\mathcal {P}(\Bbb{K},E,R)(n)=\mathcal{F}(E)(n)/\mathcal{R}(n).
\]

Son op\'erade quadratique duale est d\'efinie par
\[
\mathcal{P}^!=\mathcal{P}(\Bbb{K},E^\vee,R^\perp)
\]
o\`u $E^\vee$ est le dual de $E$ muni de l'action de $\sum_n$ duale tensoris\'ee par la repr\'esentation signature.

\medskip

\noindent 3.1. L'op\'{e}rade $\mathcal{L}ieAdm$

\smallskip

Consid\'{e}rons $R_{6}$ le $\Bbb{K}\left[ \sum_{3}\right] $-sous-module de 
$\mathcal{F} \left( E\right) (3)$
engendr\'{e} par le vecteur
\begin{eqnarray*}
u_{6} &=&x_{1}.\left( x_{2}.x_{3}\right) +x_{2}.\left( x_{3}.x_{1}\right)
+x_{3}.\left( x_{1}.x_{2}\right) -x_{2}.\left( x_{1}.x_{3}\right)
-x_{3}.\left( x_{2}.x_{1}\right) -x_{1}.\left( x_{3}.x_{2}\right)  \\
&&-\left( x_{1}.x_{2}\right) .x_{3}-\left( x_{2}.x_{3}\right) .x_{1}-\left(
x_{3}.x_{1}\right) .x_{2}+\left( x_{2}.x_{1}\right) .x_{3}+\left(
x_{3}.x_{2}\right) .x_{1}+\left( x_{1}.x_{3}\right) .x_{2}
\end{eqnarray*}
et soit $\mathcal{R}_{6}$ l'id\'{e}al de $\mathcal{F}\left( E\right) $
engendr\'{e} par $R_{6}$.

\medskip

\noindent DEFINITION 3 : L'op\'{e}rade Lie-admissible, not\'{e}e
$\mathcal{L}ieAdm$ est l'op\'{e}rade binaire quadratique 
\[
\mathcal{L}ieAdm=\mathcal{F}\left( E\right) /\mathcal{R}_{6}.
\]

\medskip

\noindent 3.2. On peut \'egalement d\'{e}finir les op\'{e}rades
binaires quadratiques associ\'{e}es \`{a} chacune des alg\`{e}bres
$G_{i}$-associatives, ($i=1,...,5$) : 

\noindent $i=1$: $\mathcal{A}ss=\mathcal{F}\left( E\right)
/\mathcal{R}_{1}$ o\`u $R_1$ est le $\Bbb{K}\left[ \sum_{3}\right]
$-sous-module engendr\'e par les vecteurs   
$$
x_{1}.\left( x_{2}.x_{3}\right)-\left(x_{1}.x_{2}\right) .x_{3}, 
$$
$\mathcal{A}ss$ est l'op\'erade des alg\`ebres associatives.

$i=2$: $\mathcal{V}inb=\mathcal{F}\left( E\right) /\mathcal{R}_{2}$ avec
${R}_{2}$ engendr\'e par \[ x_{1}.\left( x_{2}.x_{3}\right) -x_{2}.\left(
x_{1}.x_{3}\right) -\left( x_{1}.x_{2}\right) .x_{3}+\left( x_{2}.x_{1}\right)
.x_{3}, \] 
$i=3$: $\mathcal{P}reLie=\mathcal{F}\left(
E\right) /\mathcal{R}_{3}$ avec ${R}_{3}$ engendr\'e par  
\[ x_{1}.\left(
x_{2}.x_{3}\right) -x_{1}.\left(
x_{3}.x_{2}\right) -\left( x_1 .x_2 \right) .x_3 +\left( x_{1}.x_{3}\right) .x_{2}, 
\]
$i=4$: $G_{4}-\mathcal{A}ss$=$\mathcal{F}\left( E\right) /\mathcal{R}_{4}$ avec
${R}_{4}$ engendr\'e par   \[
x_{1}.\left( x_{2}.x_{3}\right) -x_{3}.\left( x_{2}.x_{1}\right)
-\left( x_{1}.x_{2}\right).x_{3} +\left( x_{3}.x_{2}\right) .x_{1},
\]
$i=5$: $G_{5}-\mathcal{A}ss$=$\mathcal{F}\left( E\right)
/\mathcal{R}_{5}$ avec ${R}_{5}$ engendr\'e par   \[
x_{1}.\left( x_{2}.x_{3}\right) +x_{2}.\left( x_{3}.x_{1}\right)
+x_{3}.\left( x_{1}.x_{2}\right) -\left( x_{1}.x_{2}\right) .x_{3}-\left(
x_{2}.x_{3}\right) .x_{1}-\left( x_{3}.x_{1}\right) .x_{2}.
\]

\bigskip

\noindent 4. {\bf Op\'{e}rades duales associ\'{e}es aux alg\`{e}bres
$G_{i}$-associatives}

\smallskip

Consid\'{e}rons le produit scalaire sur{\normalsize \ $\mathcal{F(}E)(3)$ }
pour lequel les vecteurs de base sont orthogonaux et tel que {\normalsize 
\begin{eqnarray*} &<&x_i.(x_j.x_k),x_i.(x_j.x_k)>=sgn(
\begin{tabular}{lll}
$1$ & $2$ & $3$ \\ 
$i$ & $j$ & $k$
\end{tabular}
) \\
&<&(x_i.x_j).x_k,(x_i.x_j).x_k>=-sgn(
\begin{tabular}{lll}
$1$ & $2$ & $3$ \\ 
$i$ & $j$ & $k$
\end{tabular}
)
\end{eqnarray*}
o\`u $sgn(\sigma ) =(-1)^{\varepsilon (\sigma )}$ est la signature de $\sigma$.

Soit $R_6$ le $\Bbb{K}\left[ \sum_{3}\right] $-sous-module
d\'{e}terminant l'op\'{e}rade $\mathcal{L}ieAdm$. L'orthogonal par
rapport au produit scalaire d\'{e}fini ci-dessus, $R_6 ^{\perp%
}$} , est de dimension {\normalsize $11$.} Soit {\normalsize $R_6$}$^{\prime }$
le $\Bbb{K}\left[ \sum_{3}\right] $-sous-module de 
$ \mathcal{F(}E)(3)$ engendr\'{e} par les relations  
\begin{eqnarray*}
&&(x_{\sigma (1)}x_{\sigma (2)})x_{\sigma (3)}-x_{\sigma (1)}(x_{\sigma
(2)}x_{\sigma (3)}), \\
&&(x_{\sigma (1)}x_{\sigma (2)})x_{\sigma (3)}-(x_{\sigma
(1)}x_{\sigma (3)})x_{\sigma (2)}, \\
&&(x_{\sigma (1)}x_{\sigma (2)})x_{\sigma (3)}-(x_{\sigma (2)}x_{\sigma
(1)})x_{\sigma (3)}.
\end{eqnarray*}
Alors $\dim R_6 ^{\prime }=11$ et $<u,v>=0$ pour
tout{\normalsize \ $v\in R_6^{\prime }$ , $u$ }\'{e}tant le vecteur
g\'en\'erateur de $R_6.$ Ceci implique que
$R_6 ^{\prime }\simeq R_6^{\perp }$ et par d\'efinition l'op\'erade
duale de $\mathcal{F}(E)/ {\mathcal{R}_6}= \mathcal{L}ieAdm$
not\'ee $\mathcal{L}ieAdm^{!}$ est par d\'efinition l'op\'erade
binaire quadratique $\mathcal{F(}E)/ \mathcal{R}_6^{\perp }.$ 

\medskip

\noindent PROPOSITION 1.
{\it Une alg\`ebre $A$ sur l'op\'erade $\mathcal{L}ieAdm^{!}$ est une alg\`ebre associative qui est
ab\'elienne d'ordre 3, c'est-\`a-dire v\'erifiant 
\[ abc=acb=bac
\]
pour tout $a,b,c \in A$.
  
\medskip
}
Les op\'erades duales $\mathcal{A}ss^{!}$ et $\mathcal{P}reLie^{!}$ ont
\'et\'e d\'etermin\'ees dans [4] et [3]. La premi\`ere \'etant autoduale
v\'erifie $\mathcal{A}ss^{!}=\mathcal{A}ss$ et la deuxi\`eme correspond \`a
l'op\'erade  $\mathcal{P}erm$.

\medskip

\noindent PROPOSITION 2. 
{\it Les op\'{e}rades duales de{\normalsize \ $\mathcal{V}inb,$ $G_{4}-\mathcal{A}%
ss,$ $G_{5}-\mathcal{A}ss$ }sont les op\'{e}rades quadratiques dont les
alg\`{e}bres{\normalsize \ }correspondantes sont associatives et
v\'{e}rifient respectivement : 

- pour $\mathcal{V}inb^{!}$ : $abc=bac$.

- pour $G_{4}-\mathcal{A}ss^{!}$ : $abc=cba$. 

- pour $G_{5}-\mathcal{A}ss^{!}$ : $abc=bca=cab.$ 

\medskip
}

\noindent {\it Esquisse de preuve.} $R_{2}^{\perp }$ est le {\normalsize
$\Bbb{K}\left[ \sum_{3}\right] $}-sous-module de 
$\mathcal{F(}E)(3)$ engendr\'{e} par les vecteurs 
\begin{eqnarray*}
&&\left( x_{1}.\left(
x_{2}.x_{3}\right) -\left( x_{1}.x_{2}\right) .x_{3}\right) ,\,\left(
x_{1}.\left( x_{2}.x_{3}\right) -x_{2}.\left( x_{1}.x_{3}\right) \right)
\,, \\
&&\left( x_{1}.\left( x_{2}.x_{3}\right) -\left( x_{2}.x_{1}\right)
.x_{3}\right) 
\end{eqnarray*}
 pour tout $x_{1},x_{2},x_{3}\in E$.

De m\^eme
\begin{eqnarray*}
&&R_{4}^{\perp }= \langle ( x_{1}.( x_{2}.x_{3}) -( x_{1}.x_{2}).x_{3}) ,(
x_{1}.( x_{2}.x_{3}) -x_{3}.(x_{2}.x_{1}) ) , \\ 
&& x_{1}.(x_{2}.x_{3}) -( x_{3}.x_{2}) .x_{1}) \rangle \\ 
&&R_{5}^{\perp }=\langle ( x_{1}.(x_{2}.x_{3}) -( x_{1}.x_{2}).x_{3})
,(x_{1}.( x_{2}.x_{3}) -x_{2}.( x_{3}.x_{1})) , \\ && x_{1}.( x_{2}.x_{3}) -(
x_{2}.x_{3}).x_{1}) ,( x_{1}.( x_{2}.x_{3}) -x_{3}.(x_{1}.x_{2})),( x_{1}.(
x_{2}.x_{3}) -(x_{3}.x_{1}).x_{2}) \rangle  \end{eqnarray*}
et $\dim R_{4}^{\perp }=9, \ \dim R_{5}^{\perp }=10.$ Il suffit ensuite de
remarquer que {\normalsize ($\mathcal{F(}E)/\mathcal{R} )^{!}$ }est par d\'efinition
l'op\'{e}rade binaire quadratique $\mathcal{F(}E)/\mathcal{R }^{\perp }$
d'o\`{u} le r\'{e}sultat. $\blacksquare $

\bigskip 

\noindent 5. {\bf Dualit\'e de Koszul des op\'erades provenant des alg\`ebres $G_i$-associatives.}

\smallskip

Rappelons qu'une op\'erade quadratique $\mathcal{P}$ est dite de Koszul si pour toute $\mathcal{P}$-alg\`ebre libre $F_{\mathcal{P}}(V)$ on a 
$H_{i}^{\mathcal{P}}(F_{\mathcal{P}}(V))=0, \quad i > 0$.

\medskip  

\noindent PROPOSITION 3.
{\it Les op\'erades $\mathcal{A}ss$, $\mathcal{V}inb$, $\mathcal{P}reLie$ sont de Koszul.
Les op\'erades $G_{4}-\mathcal{A}ss$ et $G_{5}-\mathcal{A}ss$ ne sont pas de Koszul.}

\medskip 

\noindent {\it Preuve.} En effet d'apr\`es [4] et [3] les op\'erades $\mathcal{A}ss$, $\mathcal{P}reLie$ sont de Koszul. Compte tenu des relations liant 
$\mathcal{P}reLie$ et $\mathcal{V}inb$, cette derni\`ere est aussi de Koszul. Pour les deux autres, nous allons montrer qu'elles ne sont pas de Koszul
en utilisant leur s\'erie de Poincar\'e et le crit\`ere d\^u \`a Ginzburg-Kapranov [4]. 
La s\'erie est d\'efinie, pour une op\'erade $\mathcal{P}$, par
\[
g _{\mathcal{P}} (x):= \sum_{i=1}^{\infty} (-1)^n dim{\mathcal{P}}(n) \frac{x^n}{n!}, 
\]

Les s\'eries de Poincar\'e d'une op\'erade de Koszul $\mathcal{P}$ et de sa duale ${\mathcal{P}} ^!$ sont reli\'ees par l'\'equation fonctionnelle 
\[
g_{\mathcal{P}}(g_{\mathcal{P} ^!}(x))=x.
\]
On a 
$$
g_{\mathcal{G}_{4}-\mathcal{A}ss} (x)= -x+x^2-\frac{3}{2}x^3+\frac{59}{4!}x^4+... \quad , \quad g_{{\mathcal{G}_{4}-\mathcal{A}ss}^{!}} (x)= -x+x^2-\frac{1}{2}x^3-\frac{1}{4}x^4+... 
$$
$$
g_{\mathcal{G}_{5}-\mathcal{A}ss} (x)= -x+x^2-\frac{10}{3!}x^3+\frac{39}{4!}x^4+... \quad , \quad g_{{\mathcal{G}_{5}-\mathcal{A}ss}^{!}} (x)= -x+x^2-\frac{1}{3}x^3+\frac{1}{12}x^4+...$$
et ces deux s\'eries ne sont pas respectivement inverses l'une de l'autre. Ces op\'erades ne sont pas de Koszul d'apr\`es [4].

\bigskip

\noindent BIBLIOGRAPHIE :

\medskip

\noindent [1] Albert A.A., On the power-associative
rings, Trans. Amer. Math. Soc. 64 (1948) 552--593. 

\noindent [2] Chapoton F., Alg\`{e}bres pr\'{e}-Lie et alg\`{e}bres
de Hopf li\'{e}es \`{a} la renormalisation, Note aux C.R.A.Sc. Paris
t.332  S\'{e}rie 1 (2001) 681-684.

\noindent[3] Chapoton F., Livernet M., Pre-Lie algebra and the rooted trees
operad, Internat. Math. Res. Notices 8 (2001) 395-408. 

\noindent[4] Ginzburg V., Kapranov M., Koszul duality for operads, Duke
Math Journal 76,1 (1994) 203-272. 

\noindent[5] Gerstenhaber M., The cohomology structure of an associative ring,
Ann of math.(2) 78 (1963) 267-288.

\noindent[6] Kobayashi S., Nomizu K., Foundations of differential geometry, Vol.I,  
John Wiley \& Sons, Inc., New York, 1963.

\noindent[7] Loday J.L., La renaissance des op\'{e}rades, S\'{e}minaire Bourbaki 
1994/95.  Ast\'{e}risque 237 (1996) 47-74.

\noindent[8] Loday J.L., Une version non commutative des alg\`ebres de Lie: les 
alg\`ebres de Leibniz, Ens. Math. 39 (1993) 269-293.

\noindent[9] May J.P., Geometry of iterated loop spaces, Lect. Notes in Math. 271, Springer-Verlag, 1972.

\noindent [10] Nijenhuis A., Sur une classe de propri\'et\'es communes a quelques types 
diff\'erents d'alg\`ebres, Enseignement math.(2) 14 (1968) 225-277. 

\noindent[11] Remm E., Structures affines sur les alg\`ebres de Lie et op\'erades Lie-admissibles, Th\`{e}se, Mulhouse 2001.

\end{document}